\documentclass[a4paper,10pt]{article}

\usepackage{amsmath,amssymb}
\usepackage[english]{babel}
\usepackage[all]{xy}

\newtheorem{lemma}{Lemma}
\newtheorem{theorem}[lemma]{Theorem}
\newtheorem{corollary}[lemma]{Corollary}
\newtheorem{proposition}[lemma]{Proposition}
\newtheorem{remark}[lemma]{Remark}
\newenvironment{proof}{\noindent{}{\emph{Proof.}}\noindent}{\hfill$\Box$ \\}

\def\mm{\mathfrak{m}}
\def\deg{\mathrm{deg}} \def\gr{\mathrm{gr}}
\def\ZZ{\mathbb{Z}}
\title{An explicit duality for quasi-homogeneous ideals}
\author{ Jean-Pierre Jouanolou \vspace {.2cm}\\
 \small Universit\'e Louis Pasteur, \\
 \small  7 rue Ren\'e Descartes, \\
 \small  67084 Strasbourg Cedex, France.\\
\small Email: {\tt jouanolo@math.u-strasbg.fr} 
} 

\begin{document}

\maketitle

\begin{abstract}
Given $r\geq n$ quasi-homogeneous polynomials in $n$ variables, the existence of a certain duality is shown and explicited in terms of generalized Morley forms. This result, that can be seen as a generalization of \cite[corollary 3.6.1.4]{Jou96} (where this duality is proved in the case $r=n$), was observed by the author at the same time. We will actually closely follow the proof of (loc.~cit.) in this paper.
\end{abstract}

\section{Notations}
Let $k$ be a non-zero unitary commutative ring. Suppose given an integer $n\geq 1$, a sequence  $(m_1,\ldots,m_n)$ of positive integers and consider the polynomial $k$-algebra 
$C:=k[X_1,\ldots,X_n]$  which is graded by setting
\begin{equation}\label{graduation}
\deg(X_i):=m_i \text{ for all } i\in \{1,\ldots,n\} \text{ and } \deg(u)=0 \text{ for all } u \in k.
\end{equation}
We will suppose moreover given an integer $r\geq n$, a sequence $(d_1,\ldots,d_r)$ of positive integers and, for all $i\in \{1,\ldots,r\}$, a (quasi-)homogeneous polynomial of degree $d_i$ 
$$f_i(X_1,\ldots,X_n):=\sum_{\substack{\alpha_1,\ldots,\alpha_n\geq 0 \\ \sum_{i=1}^n\alpha_im_i=d_i}} u_{i,\alpha}X_1^{\alpha_1}\ldots X_n^{\alpha_n} \in k[X_1,\ldots,X_n]_{d_i}.$$
In the sequel, we denote by $I$ the ideal of $C$ generated by the polynomials $f_1,\ldots,f_r$, by $\mm$ the ideal of $C$ generated by the variables $X_1,\ldots,X_n$ and by  
$B$ the quotient $C/I$. We also set $\delta:=\sum_{i=1}^rd_i -\sum_{i=1}^nm_i$.

\section{The transgression map}

Consider the Koszul complex $K^\bullet(f_1,\ldots,f_r;C)$, which is a $\ZZ$-graded complex of $C$-modules, associated to the sequence $(f_1,\ldots,f_r)$ of elements in $C$. It is of the form
$$\xymatrix{
0 \ar[r] & C(-\sum_{i=1}^rd_i) \ar@{=}[d] \ar[r]  & \cdots \ar[r] & \oplus_{i=1}^rC(-d_i) \ar@{=}[d] \ar[rr]^{\ \ \ \ \ \  (f_1,\ldots,f_r)} & & C  \ar@{=}[d]\ar[r] & 0\\
 & K^{-r} & & K^{-1} &  & K^0 & 
}$$
where, more precisely, $K^{-i}:=\bigwedge^i(K^{-1})=\bigoplus_{J\subset \{1,\ldots,r\}, |J|=i} C(-\sum_{j\in J}d_j)$. 
It gives rise to two classical spectral sequences 
$$\left\lbrace\begin{array}{lclcl}
'_1E^{pq} & := & H^q_\mm(K^p) & \Rightarrow & H^{p+q}_\mm(K^\bullet) \\
''_2E^{pq} & := & H^p_\mm(H^q(K^\bullet)) & \Rightarrow & H^{p+q}_\mm(K^\bullet). \\
\end{array}\right.$$
Since $H^i_\mm(C)=0$ if $i\neq n$, the first spectral sequence shows that, for all $p\in \ZZ$, $H^{p}_\mm(K^\bullet)$ is the cohomology module $H^{-n+p}(H^n_\mm(K^\bullet))$. Then, the second spectral sequence gives a transgression map, for all $p\in \{0,\ldots,r-n\}$,
\begin{equation}\label{trans}
H^{-n-p}(H^n_\mm(K^\bullet)) \rightarrow H^{0}_\mm(H^{-p}(K^\bullet)).
\end{equation}
In particular, taking $p=r-n$ and using the equality
$$H^{-r}(H^n_\mm(K^\bullet)) \simeq H^0(f_1,\ldots,f_r;H^n_\mm(C))(-\sum_{i=1}^rd_i),$$
we get the transgression map
$$\tau: H^0(f_1,\ldots,f_r;H^n_\mm(C))(-\sum_{i=1}^rd_i) \rightarrow H^0_\mm(H_{r-n}(K_\bullet))$$
(note that we now use the more usual homological notation for the Koszul complex: $K_p=K^{-p}$ and $H_p(K_\bullet)=H^{-p}(K^\bullet)$ for all $p\in\ZZ$) which is particularly interesting because of the
\begin{proposition}\label{isotau} If $\mathrm{depth}_I(C)= n$ then $\tau$ is an isomorphism.
\end{proposition}
\begin{proof}
If $\mathrm{depth}_I(C)= n$ then $H_i(K_\bullet)=0$ for all $i>r-n$, and then the comparison of the two spectral sequences above shows immediately that $\tau$ is an isomorphism.
\end{proof} 
\begin{remark} Observe that $\mathrm{depth}_I(C)=n$ if and only if $\mathrm{depth}_I(C)\geq n$ since $I \subset \mm$ and $\mathrm{depth}_\mm(C)=n$.  
\end{remark}

We will denote, for all $\nu \in \ZZ$, by $\tau_\nu$ the homogeneous component of degree $\nu$ of $\tau$. 
Recall that, for all $\nu\in \ZZ$ we have a canonical perfect pairing between free $k$-modules of finite type
\begin{equation}\label{perfpairing}
C_\nu \otimes_k H^n_\mm(C)_{-\nu-\sum_{i=1}^rm_i} \rightarrow H^n_\mm(C)_{-\sum_{i=1}^rm_i}\simeq k.
\end{equation}
It follows that 
$$H^0(f_1,\ldots,f_r;H^n_\mm(C))(-\sum_{i=1}^rd_i)_\nu \simeq \mathrm{Hom}_k(B_{\delta-\nu},k)=:\check{B}_{\delta-\nu}$$
(\ $\check{ }$ stands for the dual over $k$) which allows to identify $\tau_\nu$ with the $k$-modules morphism
$$\hat{\tau}_\nu:\check{B}_{\delta-\nu} \rightarrow H^0_\mm(H_{r-n}(K_\bullet))_\nu.$$
Observe that the direct sum $\bigoplus_{\nu \in \ZZ} \check{B}_{\delta-\nu}$ has a natural structure of $B$-module and so the $k$-linear map
$$\hat{\tau}:=\bigoplus_{\nu\in\ZZ}\hat{\tau}_\nu : \bigoplus_{_\nu \in \ZZ} \check{B}_{\delta-\nu} \rightarrow H^0_\mm(H_{r-n}(K_\bullet))$$
is a morphism of graded $B$-modules. It is clear that $\hat{\tau}$ is an isomorphism if $\tau$ is itself an isomorphism, and for instance if $\mathrm{depth}_I(C)\geq n$ by proposition \ref{isotau}. In the rest of this note we will give an explicit description of the map $\hat{\tau}$ in this case.

\section{Generalized Morley forms}

Introducing new indeterminates $Y_1,\ldots,Y_n$, we  identify the ring $C\otimes_k C$ with the polynomial ring $k[\underline{X},\underline{Y}]$ (we shortcut sequences: for instance $\underline{X}$ stands for the sequence $(X_1,\ldots,X_n)$) which is canonically graded via the tensor product: $\deg(X_i)=\deg(Y_i)=m_i$ for all $i=1,\ldots,r$. 

 In $C\otimes_k C$, for all $i\in\{1,\ldots,r\}$  we choose a decomposition
\begin{eqnarray}\label{decomp} 
f_i(X_1,\ldots,X_n)-f_i(Y_1,\ldots,Y_n) & = & f_i\otimes_k 1-1\otimes_k f_i \\ \nonumber
& = & \sum_{j=1}^n 
(X_j-Y_j)g_{i,j}(X_1,\ldots,X_n,Y_1,\ldots,Y_n).
\end{eqnarray}
Let $e_1,\ldots,e_r$ be the canonical basis of $\oplus_{i=1}^rC\otimes_kC(-d_i)$ with $\deg(e_i)=d_i$ for all $i\in\{1,\ldots,r\}$ and consider 
$$\Delta:=\sum_{\substack{\sigma \in \mathfrak{S}_r \text{ such that}\\ \sigma(1)<\cdots<\sigma(n), \\ \sigma(n+1)<\cdots<\sigma(r)}} \epsilon(\sigma) 
\left|\begin{array}{cccc}
g_{\sigma(1),1} & g_{\sigma(1),2} & \cdots & g_{\sigma(1),n} \\
g_{\sigma(2),1} & g_{\sigma(2),2} & \cdots & g_{\sigma(2),n} \\
\vdots & \vdots & & \vdots \\
g_{\sigma(n),1} & g_{\sigma(n),2} & \cdots & g_{\sigma(n),n} \\
      \end{array}
\right|
e_{\sigma(n+1)}\wedge\cdots\wedge e_{\sigma(r)}$$
where $\mathfrak{S}_r$ denotes the set of all the permutations of $r$ elements and $\epsilon(\sigma)$ the signature of such a permutation $\sigma \in \mathfrak{S}_r$.

\begin{lemma}[{\cite[2.14.2]{Jou80}}] The element $\Delta$ is a cycle of the Koszul complex associated to the sequence $(f_1\otimes_k 1-1\otimes_k f_1,\ldots,f_r\otimes_k 1-1\otimes_k f_r)$ in $C\otimes_kC$.
\end{lemma}
We denote by $\Delta^\flat$ the class of $\Delta$ in the homology group
$$H_{r-n}(f_1(\underline{X})-f_1(\underline{Y}),\ldots,f_r(\underline{X})-f_r(\underline{Y});C\otimes_kC)_\delta.$$Note that, as a consequence of the so-called Wiebe lemma (see for instance \cite[3.8.1.7]{Jou95}), $\Delta^\flat$ does not depend on the choice of the decompositions \eqref{decomp} since the sequence  $(X_1-Y_1,\ldots,X_n-Y_n)$ is regular in $C\otimes_kC$.

\begin{lemma}\label{DP} For all $P \in C$ we have $(P\otimes_k 1-1\otimes_k P)\Delta^\flat=0$, or in other words
$P(\underline{X})\Delta^\flat=P(\underline{Y})\Delta^\flat$. 
\end{lemma}
\begin{proof} First, it is clear that for all $i=1,\ldots,n$ we have $(X_i-Y_i)\Delta^\flat=0$. Indeed, each determinant fitting in the definition of $\Delta$ becomes an element of the ideal generated by the polynomials $f_j(\underline{X})-f_j(\underline{Y})$, $j=1,\ldots,r$, after multiplication by $(X_i-Y_i)$. The proof then follows from the equality
$$P(\underline{X})-P(\underline{Y})  =  \sum_{i=1}^{n} P(Y_1,\ldots,Y_{i-1},X_i,\ldots,X_n)-P(Y_1,\ldots,Y_{i},X_{i+1},\ldots,X_n)$$
where each term in the above sum is divisible by at least one of the elements $(X_i-Y_i)$,  $i\in\{1,\ldots,n\}$. 
\end{proof}

The canonical projection $C\otimes _k C \rightarrow C \otimes _k B$ induces a map
\begin{equation}\label{CB}
\xymatrix{
 H_{r-n}(f_1\otimes_k 1-1\otimes_k f_1,\ldots,f_r\otimes_k 1 -1\otimes_k f_r;C\otimes_kC) \ar[d] \\
H_{r-n}(f_1\otimes_k 1,\ldots,f_r\otimes_k 1 ;C\otimes_kB) 
}
\end{equation}
(note that $1\otimes_k f_i=0$ in $C\otimes_k B$ for all $i=1,\ldots,r$) which sends $\Delta^\flat$ to an element, that we will denote $\nabla$, of degree $\delta$ in $H_{r-n}(f_1\otimes_k 1,\ldots,f_r\otimes_k 1 ;C\otimes_kB).$

Observe that, for all $q\in\ZZ$ the $C$-module 
$H_{r-n}(f_1\otimes_k 1,\ldots,f_r\otimes_k 1;C\otimes_k B_q)$ is $\ZZ$-graded via the grading of $C$, so we deduce that the $B\otimes_k B$-module 
$H_{r-n}(f_1\otimes_k 1,\ldots,f_r\otimes_k 1 ;C\otimes_kB)$ is bi-graded; for all $p,q \in \ZZ\times \ZZ$ we have
$$H_{r-n}(f_1\otimes_k 1,\ldots,f_r\otimes_k 1 ;C\otimes_kB)_{p,q} := H_{r-n}(f_1\otimes_k 1,\ldots,f_r\otimes_k 1;C\otimes_k B_q)_p.$$
We can thus decompose $\nabla$ with respect to this bi-graduation and  we define $\nabla=\sum_{(p,q) \in \ZZ^2}\nabla_{p,q}$ with
$$\nabla_{p,q} \in H_{r-n}(f_1\otimes_k 1,\ldots,f_r\otimes_k 1 ;C\otimes_k B_q)_{p}.$$

\begin{lemma} For all couple $(p,q)\in \ZZ^2$ we have
$$\nabla_{p,q} \in H^0_{\mm\otimes_kB+B\otimes _k \mm}(H_{r-n}(f_1\otimes_k 1,\ldots,f_r\otimes_k 1 ;C\otimes_kB))_{p,q}.$$
\end{lemma}
\begin{proof} This lemma follows from lemma \ref{DP}; for all $j\in \{1,\ldots,n\}$ we have the equality $(X_j\otimes_k 1 -1\otimes_k X_j)\nabla=0$ which gives, by looking at the homogeneous components, 
$$(X_j\otimes_k 1)\nabla_{p,q} = (1\otimes_k X_j)\nabla_{p+1,q-1}$$
for all $(p,q)\in \ZZ^2$ such that $p+q=\delta.$ By successive iterations we obtain
$$(X_j\otimes_k 1)^{q+1}\nabla_{p,q}=(1\otimes_k X_j)^{q+1}\nabla_{p+q-1,-1}=0$$
which shows that $(\mm\otimes_k B)^{nq+1}\nabla_{p,q}=0$. Exactly in the same way we obtain $(B\otimes_k \mm)^{np+1}\nabla_{p,q}=0$.
\end{proof}

Finally, let us emphasize that $\nabla_{\delta,0}$ has a simple description. For all $i\in\{1,\ldots,r\}$ we  choose a decomposition 
\begin{equation}\label{decompX}
f_i(X_1,\ldots,X_n)=\sum_{j=1}^n X_jf_{i,j}(X_1,\ldots,X_n) \in C
\end{equation}
and similarly to what we did above, we consider 
$$\Lambda:=\sum_{\substack{\sigma \in \mathfrak{S}_r \text{ such that}\\ \sigma(1)<\cdots<\sigma(n), \\ \sigma(n+1)<\cdots<\sigma(r)}} \epsilon(\sigma) 
\left|\begin{array}{cccc}
f_{\sigma(1),1} & f_{\sigma(1),2} & \cdots & f_{\sigma(1),n} \\
f_{\sigma(2),1} & f_{\sigma(2),2} & \cdots & f_{\sigma(2),n} \\
\vdots & \vdots & & \vdots \\
f_{\sigma(n),1} & f_{\sigma(n),2} & \cdots & f_{\sigma(n),n} \\
      \end{array}
\right|
e_{\sigma(n+1)}\wedge\cdots\wedge e_{\sigma(r)}.$$
It is, as $\Delta$, a cycle of the Koszul complex $K_\bullet(f_1,\ldots,f_r;C)$. We denote  
$\bar{\Lambda}$ its class in $H_{r-n}(K_\bullet(f_1,\ldots,f_r;C))_\delta$, class which is independent, by the Wiebe lemma, of the choice of the decompositions \eqref{decompX} since the sequence $(X_1,\ldots,X_n)$ is regular in $C$.
\begin{lemma}\label{d0} We have $\nabla_{\delta,0}=\bar{\Lambda}$ in $H_{r-n}(K_\bullet(f_1,\ldots,f_r;C))_\delta$.
\end{lemma}
\begin{proof} Indeed, $\nabla_{\delta,0}$ is the image of $\nabla^\flat$ via the map
$$C\otimes_kC=k[\underline{X},\underline{Y}] \rightarrow C=k[\underline{X}] : P(\underline{X},\underline{Y}) \mapsto P(\underline{X},0)$$
and this shows immediately the claimed equality.
\end{proof}

\section{The explicit duality}

Suppose given $\nu \in \ZZ$ and $u\in \check{B}_{\delta-\nu}=\mathrm{Hom}_k(B_{\delta-\nu},k)$. The canonical morphism 
$\mathrm{id}_C\otimes_k u:C\otimes_k B_{\delta-\nu} \rightarrow C\otimes_k k\simeq C$ induces a map 
$$H_{r-n}(f_1,\ldots,f_r;C\otimes_k B_{\delta-\nu})_\nu \xrightarrow{H_{r-n}(f_1,\ldots,f_r;\mathrm{id}_C\otimes_k u)} 
H_{r-n}(f_1,\ldots,f_r;C)_\nu$$
which sends $\nabla_{\nu,\delta-\nu}$ to an element that we will denote $\nabla^{(u)}_{\nu,\delta-\nu}$. Therefore, to any $u\in \check{B}_{\delta-\nu}$ we can associate an element in $H_{r-n}(f_1,\ldots,f_r;C)_\nu$. Denoting $D^{\gr}_k(B)$ the graded $B$-module of graded morphisms from $B$ to $k$, that is to say
	$$D^{\gr}_k(B):=\mathrm{Hom}_k^{\gr}(B,k)=\bigoplus_{\nu \in \ZZ} \mathrm{Hom}_k(B,k)_\nu=\bigoplus_{\nu \in \ZZ} \mathrm{Hom}_k(B_{-\nu},k)=\bigoplus_{\nu \in \ZZ} \check{B}_{-\nu}$$
we obtain a map
\begin{eqnarray}\label{omega}
\omega : D^{\gr}_k(B)(-\delta) & \rightarrow & H_{r-n}(f_1,\ldots,f_r;C) 
\end{eqnarray}
and we have the

\begin{proposition}\label{Blinear} The map $\omega$ is a graded morphism (i.e.~of degree 0) of graded $B$-modules whose image 
is contained in $H^0_\mm(H_{r-n}(f_1,\ldots,f_r;C))$. 
\end{proposition}
\begin{proof} 
Let us choose a couple $(q,\nu) \in \ZZ^2$ and pick up $b\in B_q$ and $u \in D^{\gr}_k(B)(-\delta)_\nu=\check{B}_{\delta-\nu}$.  To prove the $B$-linearity of $\omega$ we have to prove that $\omega(bu)=b\omega(u)$.

On the one hand, $bu \in \check{B}_{\delta-\nu-q}$ so $\omega(bu) \in H_{r-n}(\underline{f};C)_{\nu+q}$ is, by definition, the image of $\nabla_{\nu+q,\delta-\nu-q}$ by the map
$$H_{r-n}(\underline{f};C\otimes_k B_{\delta-\nu-q})_{\nu+q} \rightarrow H_{r-n}(\underline{f};C)_{\nu+q}$$
induced by $C\otimes_kB_{\delta-\nu-q}\rightarrow C: c\otimes_k x \mapsto cu(bx)$, which is also the image of
$(1\otimes_k b)\nabla_{\nu+q,\delta-\nu-q} \in H_{r-n}(\underline{f};C\otimes_k B_{\delta-\nu})_{\nu+q}$ by the map
$$H_{r-n}(\underline{f};C\otimes_k B_{\delta-\nu})_{\nu+q} \rightarrow H_{r-n}(\underline{f};C)_{\nu+q}$$
induced by $C\otimes_kB_{\delta-\nu}\rightarrow C: c\otimes_k y \mapsto cu(y)$.

On the other hand, $b\omega(u)$ is the image of $\nabla_{\nu,\delta-\nu}$ by the map
$$H_{r-n}(\underline{f};C\otimes_k B_{\delta-\nu})_{\nu} \rightarrow H_{r-n}(\underline{f};C)_{\nu}$$
induced by $C\otimes_kB_{\delta-\nu}\rightarrow C: c\otimes_k x \mapsto c_1cu(x)$ where $c_1 \in C$ is such that $c_1=b$ in $B=C/I$. It follows that $b\omega(u)$ is the image of $(b\otimes_k 1)\nabla_{\nu,\delta-\nu}$ by the map
$$H_{r-n}(\underline{f};C\otimes_k B_{\delta-\nu})_{\nu+q} \rightarrow H_{r-n}(\underline{f};C)_{\nu+q}$$
induced by $C\otimes_kB_{\delta-\nu}\rightarrow C: c\otimes_k x \mapsto cu(x)$.

Now, by lemma \ref{DP}, we know that $(1\otimes_k b-b\otimes_k 1)\nabla=0$ which implies, looking at the homogeneous component of degree $(\nu+q,\delta-\nu)$, that 
$$b\omega(u)=(b\otimes_k 1) \nabla_{\nu,\delta-\nu}=(1\otimes_k b)\nabla_{\nu+q,\delta-\nu-q}=\omega(bu).$$

Finally, we have $D^{\gr}_k(B)(-\delta)_\nu=0$ for all $\nu>\delta$ so the $B$-linearity of $\omega$ implies that $B_q\mathrm{Im}(\omega)=0$ for all sufficiently large integer $q$, which is equivalent to $\mm^p\mathrm{Im}(\omega)=0$ in $H_{r-n}(\underline{f};C)$ for all sufficiently large integer $p$.
\end{proof}

According to the above proposition \ref{Blinear}, and abusing notation, from now on we will assume that $\omega$ denotes the map \eqref{omega} co-restricted to $H^0_\mm(H_{r-n}(f_1,\ldots,f_r;C)).$  We also define $\omega_\nu$ as the graded component of degree $\nu$ of $\omega$:
\begin{eqnarray*}
\omega_\nu : D^{\gr}_k(B)(-\delta)_\nu=\check{B}_{\delta-\nu} & \rightarrow & H^0_\mm(H_{r-n}(f_1,\ldots,f_r;C))_\nu \\
 u & \mapsto & \nabla^{(u)}_{\nu,\delta-\nu}.
\end{eqnarray*}
We are now ready to state the main result of this note.

\begin{theorem} If $\mathrm{depth}_I(C)= n$ then $\hat{\tau}=\omega$.
\end{theorem}
\begin{proof} We will prove that $\hat{\tau}_\nu=\omega_\nu$ for all $\nu \in \ZZ$. Recall that under the hypothesis $\mathrm{depth}_I(C)\geq n$ the map $\tau$ , and hence $\hat{\tau}$, become an isomorphism.

First, since $H_{r-n}(f_1,\ldots,f_r;C)$ is a sub-quotient of 
$$\bigwedge^{r-n}\left(\bigoplus_{i=1}^rC(-d_i)\right)=\bigoplus_{1\leq i_1< i_2 < \cdots < i_{r-n}\leq r} C(-d_{i_1}-d_{i_2}-\cdots -d_{i_{r-n}}),$$ 
we deduce that $H_{r-n}(f_1,\ldots,f_r;C)_\nu=0$ for all $v<0$ (note that the extreme case is obtained when $r=n$). It follows that $\omega_\nu$ and $\hat{\tau}_\nu$ are both the zero map if $\nu<0$. 

If $\nu > \delta$ then $H^0_\mm(C)(-\sum_{i=1}^rd_i)_\nu=0$. Since by hypothesis, $\hat{\tau}$ is an isomorphism 
we deduce that $$H^0_\mm(H_{r-n}(f_1,\ldots,f_r;C))_\nu=0 $$ 
and hence that $\omega_\nu$ and $\hat{\tau}_\nu$ are again both the zero map.

We now prove that $\omega_\delta=\hat{\tau}_\delta$. By definition, 
$$\omega_\delta : \check{B}_0\simeq k \rightarrow H^0_\mm(H_{r-n}(\underline{f};C))_\delta$$
sends any $\lambda \in k$ to $\lambda \nabla_{\delta,0}=\lambda \bar{\Lambda}$ (see lemma \ref{d0}), so it is completely determined by the formula $\omega_\delta(1)=\bar{\Lambda}$. To explicit the map $\hat{\tau}_\nu$ we will use the functoriality property of $\tau$ (and hence of $\hat{\tau}$); in this order, we will specify the sequence $\underline{f}:=(f_1,\ldots,f_r)$ or $\underline{X}:=(X_1,\ldots,X_n)$ in $C$ under consideration with the obvious  notation $\tau(\underline{f})$ or $\tau(\underline{X})$. The decomposition \eqref{decompX} gives a graded morphism of $C$-modules (recall $r\geq n$)
$$ \oplus_{i=1}^r C(-d_i) \xrightarrow{M} \oplus_{i=1}^n C(-m_i)$$
which can be lifted to a graded morphism of complexes
$$K_\bullet(f_1,\ldots,f_r;C) \xrightarrow{\wedge^\bullet{M}} K_\bullet(X_1,\ldots,X_n;C).$$
Note that $\wedge^p(M)$ is the zero map for all $p > n$.  
Using the self-duality property of the Koszul complexes, we obtain by duality a graded morphism of graded complexes
$$K(X_1,\ldots,X_n;C)(-\delta) \rightarrow K(f_1,\ldots,f_r;C)$$
which is of the form, denoting $K_1:=\oplus_{i=1}^r C(-d_i)$,
$$\xymatrix{
C(-\sum_{i=1}^r d_i) \ar[r] \ar[d]^{\mathrm{id}}& \cdots \ar[r]^{\underline{X}} &  C(-\delta) \ar[r] \ar[d]^{``\Lambda``}& 0 \ar[r] \ar[d]^{0}& \cdots \ar[r] & 0 \ar[d]^{0}\\
C(-\sum_{i=1}^r d_i) \ar[r] & \cdots \ar[r] &  \wedge^{r-n}K_1 \ar[r] & \wedge^{r-n-1}K_1 \ar[r] & \cdots \ar[r]^{\underline{f}} & C \\
}$$ 
By functoriality of the transgression map $\tau$ for morphisms of complexes, we obtain the commutative diagram
$$\xymatrix{ H^0(X_1,\ldots,X_n;H^n_\mm(C))(-\sum_{i=1}^r d_i) \ar[r]^{\hspace{1cm}\tau(\underline{X})(-\delta)} \ar[d]^{\mathrm{id}} & H^0_\mm(C/\mm)(-\delta)=k(-\delta) \ar[d]^{1 \mapsto \Lambda} \\
H^0(f_1,\ldots,f_r;H^n_\mm(C))(-\sum_{i=1}^r d_i) \ar[r]^{\hspace{.6cm} \tau(\underline{f})} & H^0_\mm(H_{r-n}(f_1,\ldots,f_r;C))
}$$
which yields in degree $\delta$ the commutative diagram
$$\xymatrix{ k=H^n_\mm(C)_{-\sum_{i=1}^n m_i} \ar[r]^{\hspace{1cm}\hat{\tau}_0(\underline{X})} \ar[d]^{\mathrm{id}} & k \ar[d]^{1 \mapsto \Lambda} \\
k=H^n_\mm(C)_{-\sum_{i=1}^n m_i} \ar[r]^{\hat{\tau}_\delta(\underline{f})\hspace{.45cm}} & H^0_\mm(H_{r-n}(f_1,\ldots,f_r;C))_\delta
}$$
Since the map $\hat{\tau}_0(\underline{X})$ is the identity \cite[2.6.4.6]{Jou95}, we deduce that for all $\lambda \in k$ we have
$$\hat{\tau}_\delta(\underline{f})(\lambda)=\lambda \bar{\Lambda} \in H^0_\mm(H_{r-n}(f_1,\ldots,f_r;C))_\delta$$
and hence that $\hat{\tau}_\delta=\omega_\delta$.

Finally, assume that $0\leq \nu <\delta $. By $B$-linearity of $\hat{\tau}$ and $\omega$ (see proposition \ref{Blinear}), for all $u \in D^{\gr}_k(B)(-\delta)_\nu=\check{B}_{\delta-\nu}$ and for all $ b \in B_{\delta-\nu}$ we have
$$b(\hat{\tau}_\nu(u)-\omega_\nu(u))=\hat{\tau}_\delta(bu)-\omega_\delta(bu)=0 \in H^0_\mm(H_{r-n}(f_1,\ldots,f_r;C))_\delta$$
with $\hat{\tau}_\nu(u)-\omega_\nu(u) \in H^0_\mm(H_{r-n}(f_1,\ldots,f_r;C))_\nu$.
Since $H_{r-n}(\underline{f};C)$ is a $B$-module, we have, for all $\nu\in \ZZ$ a canonical $k$-linear pairing
$$B_{\delta-\nu} \otimes_k H_{r-n}(\underline{f};C)_\nu \rightarrow H_{r-n}(\underline{f};C)_\delta.$$
By hypothesis, $\tau$ is an isomorphism and therefore we have the commutative diagram
$$\xymatrix{
B_{\delta-\nu} \otimes_k H^0(\underline{f};H^n_\mm(C))_{\nu-\sum_{i=1}^r d_i} \ar[r] \ar[d]^{\mathrm{id}\otimes_k \tau_\nu}_{\wr}&
H^0(\underline{f};H^n_\mm(C))_{-\sum_{i=1}^n m_i} \ar[d]^{\mathrm{id}\otimes_k \tau_\delta}_{\wr} \\
B_{\delta-\nu} \otimes_k H^0_\mm(H_{r-n}(\underline{f};C))_{\nu} \ar[r] &
H^0_\mm(H_{r-n}(\underline{f};C))_{\delta}}$$
where both vertical arrows are isomorphisms. Now, the top row being a non-degenerated pairing by 
\eqref{perfpairing}, we deduce that the bottom row is also a non-degene\-ra\-ted pairing and hence that $\hat{\tau}_\nu=\omega_\nu$, as claimed.
\end{proof}

\begin{corollary} If $\mathrm{depth}_I(C)= n$ then $\omega$ is an isomorphism of $B$-modules. 
\end{corollary}

\section*{Acknowledgment}

 I thank heartly Laurent Bus\'e who proposed me to write down my notes and to translate them into english. He did a good job, but refused to be a co-author of the paper.


\end{document}